\input amstex.tex
\loadeusm
\loadeusb
\documentstyle{amsppt}

\topmatter

\title
Countable dense homogeneity of definable spaces
\endtitle

\author
Michael Hru\v s\' ak and
Beatriz Zamora Avil\'{e}s
\endauthor

\thanks {\it Key words and phrases}\rm : Countable dense homogeneous,
Borel, Baire
\smallskip
\smallskip
The first author's research  was supported partially by
grant GA\v CR 201/03/0933 and by a PAPIIT grant IN108802-2 and CONACYT grant 40057-F.
\endthanks

\date June 13, 2003
\enddate

\abstract We investigate which definable separable metric spaces are countable dense homogeneous (CDH).
We prove that a Borel CDH space is completely metrizable and give a complete list of zero-dimensional
Borel CDH spaces. We also show that for a Borel $X\subseteq 2^\omega$ the following are equivalent:
(1) $X$ is $G_\delta$ in $2^\omega$, (2) $X^\omega$ is CDH and (3) $X^\omega$ is homeomorphic
to $2^\omega$ or to $\omega^\omega$. Assuming the Axiom of Projective Determinacy the results
extend to all projective sets and under the Axiom of Determinacy to all separable metric spaces. 
In particular, modulo large cardinal assumption  it is relatively consistent with ZF that all
CDH separable metric spaces are completely metrizable.
We also answer a question of Stepr$\bar{\text{a}}$ns and Zhou
by showing that $\frak p= \min \{\kappa: 2^\kappa$ is not CDH$\}$.

\endabstract
\endtopmatter

\head 0. Introduction \endhead

\bigskip

A separable topological space $X$ is {\it countable dense homogeneous (CDH)} if given any
two countable dense subsets $D, D' \subseteq X$ there is a homeomorphism $h$ of $X$ such that
$h[D]=D'$. The first result in this area is due to Cantor, who, in effect, showed that
the reals are CDH. Fr\'{e}chet [Fr] and Brower [Br], independently, proved that the same is true for
the $n$-dimensional Euclidean space $\Bbb R^n$. In 1962, Fort [Fo] proved that the
Hilbert cube is also CDH. \par
Systematic study of CDH spaces was initiated by Bennett [Be] in 1972. Since then a number
of papers were published on the topic, most of which are mentioned in the references. 
The focus remained on separable metric spaces. Under some set-theoretic assumptions like the 
Continuum Hypothesis or Martin's Axiom a variety of examples of countable dense homogeneous metric 
spaces were constructed: Assuming CH Fitzpatrick and Zhou constructed a CDH Bernstein subset of 
$\Bbb R^n$ and a CDH subset of $\Bbb R$ which is meager in itself; Baldwin and Beaudoin constructed
Bernstein subset of $\Bbb R$ under  Martin's Axiom for countable partial orders.\par

In this paper we are concerned mostly with countable dense homogeneity of definable separable metric spaces.
Our principal result
states that every analytic CDH space is completely Baire. We use it to 
give a complete list of zero-dimensional
Borel CDH spaces and to show that for a Borel $X\subseteq 2^\omega$ the following are equivalent:
(1) $X$ is $G_\delta$ in $2^\omega$, (2) $X^\omega$ is CDH and (3) $X^\omega$ is homeomorphic
to $2^\omega$ or to $\omega^\omega$. These provide partial answers to the following problems of [FZ3]:
\smallskip
{\bf 387.} For which 0-dimensional subsets of $\Bbb R$ is $X^\omega$ homogeneous? CDH?
\smallskip
\flushpar
and
\smallskip
{\bf 389.} Does there exist a CDH metric space that is not completely metrizable?

\bigskip

\head 1. Descriptive set theory \endhead

\bigskip

In this section we review some of the classical results of descriptive set theory. For proofs and
further reference consult e.g. [Ke]. Recall that a separable completely metrizable space is called
a {\it Polish} space. We call a separable metric space {\it Borel}, if it
is Borel in its completion. A separable metric space is {\it analytic} if it is a continuous
image if the Baire space $\omega^\omega$. A space is  {\it co-analytic} if it is a
complement of an analytic subspace of some Polish space. Recall that a space is Borel
if and only if it is both analytic and co-analytic. This is an old result of Souslin as is the following:
\bigskip
\proclaim{Theorem 1.1} Every analytic space contains a homeomorphic copy of $2^\omega$.
\endproclaim
\bigskip
Recall that a subset $A$  of a Polish space $X$ is said to have the {\it Baire property} if there
is an open set $U \subseteq X$ such that the symmetric difference $A\triangle U$ is meager in $X$.
\bigskip
\proclaim{Theorem 1.2}Every analytic subspace of a Polish space has the Baire property.
\endproclaim
\bigskip
A topological space $X$ is {\it Baire} if the complement of every meager subset of $X$ is dense in $X$.
Note that being Baire and having the Baire property are quite different notions. We will use the following
corollary of Theorem 1.2:
\bigskip
\proclaim{Theorem 1.3} Every analytic Baire space has a dense completely metrizable subspace.
\endproclaim

\demo{Proof} Let $X$ be an analytic Baire space and let $\bar X$ be its completion. By Theorem 1.2
there is an open set $U\subseteq \bar X$ such that $X\triangle U$ is meager in $\bar X$. That is
$X\triangle U =\bigcup_{n\in\omega} F_n$, where each $F_n$ is nowhere dense in $\bar X$. Note that 
$U$ is a dense open subset of $\bar X$.
Let $G=U\setminus \bigcup _{n\in\omega}\bar F_n$. Then $G$ is completely metrizable as it is $G_\delta$
in $\bar X$, and $G$ is a dense subset of $X$ as $X$ is Baire.\hfill\qed
\enddemo
\bigskip
A topological space $X$ is {\it completely Baire} if all of its closed subspaces are Baire. The following theorem
is due to Hurewicz (see [Ke]).

\bigskip
\proclaim{Theorem 1.4} Every co-analytic completely Baire space is completely metrizable.
\endproclaim
\bigskip

Under the Axiom of Projective Determinacy (PD) all of the above theorems hold for all projective sets.
Similarly under the Axiom of Determinacy (AD) they hold for all separable metric spaces. For proof 
of the analogues of Theorems 1.1 and 1.2 (and hence also 1.3) in this context see e.g Theorem 27.9
of [Ka]. The fact that the variants of the Theorem 1.4 hold follows from the proof of Theorem 4
of [KLW].

\bigskip

The following characterization of zero-dimensional Polish spaces can be found in [Ke] and [vM2] 
\bigskip
\proclaim{Theorem 1.5} (i) Every zero-dimensional separable compact completely metrizable space without isolated points
is homeomorphic to $2^\omega$.\newline
(ii) Every zero-dimensional separable locally compact non-compact completely me{\-}trizable space without isolated points
is homeomorphic to $2^\omega\setminus \{0\}$.\newline
(i) Every zero-dimensional separable completely metrizable space without isolated points in which all compact sets
are nowhere dense is homeomorphic to $2^\omega$.

\endproclaim

\bigskip

\head 2. Analytic CDH spaces
\endhead

\bigskip

In the article {\it Some Open Problems in Densely homogeneous spaces} of the {\it Open problems in topology} 
Fitzpatrick and Zhou ask (Question {389.}) 
whether there is a 
CDH metric space which is not completely metrizable. We answer this question in the negative for Borel spaces. 
The following simple lemma ([FZ2]) will be used many times in what follows.

\bigskip
\proclaim{Lemma 2.1} A separable metric space $X$ without isolated points is meager in itself if and only if  
there is a countable dense $D\subseteq X$ which is $G_\delta$ in $X$.
\endproclaim

\demo{Proof} The reverse implication is obvious. For the forward implication let $X=\bigcup_{n\in\omega} F_n$, where
each $F_n$ is a closed nowhere dense subset of $X$. Enumerate a
basis for the topology of $X$ as $\{U_n: n\in\omega\}$ and recursively pick $x_n\in U_n \setminus  \bigcup_{m\leq n} F_m$.
Set $D=\{x_n: n\in\omega\}$. $D$ is obviously a countable dense subset of $X$. To see that it is $G_\delta$ in $X$
note that $D$ intersects each $F_n$ in a finite set, hence $X\setminus D=\bigcup_{n\in\omega} (F_n\setminus D)$ is $F_\sigma$
in $X$.\hfill\qed
\enddemo
\bigskip

Next we prove a decomposition lemma for CDH spaces.

\bigskip
\proclaim{Lemma 2.2} Every CDH space $X$ can be written as a disjoint topological sum $X= I \oplus L\oplus R$,
where $I$ is the set of isolated points in $X$, $L$ is locally compact without isolated points and
$R$ has the property that every compact subset or $R$ is nowhere dense in $R$.
\endproclaim

\demo{Proof} First we show that the set $I$ of all isolated points of $X$ is clopen in $X$. Note that $I$ is countable as
$X$ is separable. If $I$ is not closed, pick $x\in \bar I \setminus I$
and a set $C\subseteq X\setminus \bar I$ countable dense in $X\setminus \bar I$. Let $D_0= I\cup C$ and
$D_1=D_0\cup\{x\}$. The sets $D_0$ and $D_1$ are then countable dense subsets of $X$ and we reach a contradiction by noting that
there is no homeomorphism of $X$ sending  $D_1$ to $D_0$, for
$x$ is not isolated but every neighborhood of $x$ contains an isolated point, whereas
all points in $D_0$ are either isolated or have a neighborhood which does not contain any isolated points.
\smallskip
Let $Y=X\setminus I$. Consider
\smallskip
\flushpar
$L=\{x\in Y: \exists U\subseteq Y\text{ locally compact neighborhood of }x\}\text{ and} $\newline
$R=\{x\in Y: \exists\ U\subseteq Y\text{ neighborhood of } x,\text{ s. t. } \forall\  K\subseteq U \text{ compact } 
\text{ int}(K)=\emptyset\}.$
\smallskip
\flushpar
Obviously $L$ and $R$ are disjoint open subsets of $Y$.  To finish the proof it suffices to show that
$Y=L\cup R$. First note that $L\cup R$ is dense in $Y$, as if $x\in Y\smallsetminus (L\cup R)$ then 
i.p. $x\in Y\smallsetminus R$, which implies that for every $ U\subseteq X$ neighborhood of $x$, there is a $K\subseteq U$ compact 
such that $\text{int}(K)\neq\emptyset$, hence $x\in \overline {L}$.\par

Now, suppose that $Y\smallsetminus (L\cup R)\neq\emptyset$. Pick $x\in Y\smallsetminus (L\cup R)$ and a 
countable dense $D_0\subseteq L\cup R$ and let $D_1=D_0\cup \{x\}$. Again, $D_0$ and $D_1$ are clearly countable in $X$ 
and there is no homeomorphism $h$
 of $X$ sending  $D_1$ to $D_0$ as then $h(x)\in L$ or $h(x)\in R$ but $x\not \in L\cup R$.\hfill\qed
\enddemo

\bigskip
\proclaim{Theorem 2.3} Every analytic CDH space $X$ is completely Baire.
\endproclaim

\demo{Proof} By Lemma 2.2 we can assume that $X$ has no isolated points.
\smallskip
{\bf Claim 1.} Every open subset of $X$ is uncountable.
\smallskip
  
Assume not, that is $V=\bigcup\{U: U$ is a countable open subset of $X\}$ is not empty. Then
$V$ is itself a countable open set. Choose $C$ a countable dense subset of $X\setminus V$ and $x\in V$.
Let $D_0= C\cup V$ and $D_1=C\cup V\setminus \{x\}$. The sets $D_0$ and $D_1$ are then countable dense subsets of $X$.
As $X$ is CDH there is a homeomorphism $h$ of $X$ such that $h[D_1]=D_0$. Then, however, $h(x)\not \in V$ and,
unlike $x$, $h(x)$ does not have a countable neighborhood which contradicts the fact that $h$ is a homeomorphism.
\smallskip
{\bf Claim 2.} $X$ is Baire.
\smallskip
  
 Suppose it is not the case. That means that there is an open set $U\subseteq X$ which is meager in itself. By Lemma 2.1
there is a $C\subseteq U$ countable dense in $U$ which is $G_{\delta}$ in $U$. Let $D_0$ be a countable dense subset of
$X$ such that $D_0\cap U= C$. \par
Let $\{U_{n}:n\in\omega\}$ be an enumeration of some countable basis for the topology on $X$. By Claim 1, each $U_n$ is
uncountable, as every open subset of an analytic space is itself analytic, by Theorem 1.1, each $U_n$ contains
a subset $F_n$ homeomorphic to $2^\omega$. Choose, for every $n\in\omega$, a countable $C_n\subseteq F_n$ dense in $F_n$ and 
set $D_1=\bigcup_{n\in\omega} C_n$. The set $D_1$ is then a countable dense subset of $X$. \par
Note that $D_1 \cap V$ is not  $G_{\delta}$ in $V$ for any open set $V\subseteq X$.
To see this let $V$ be an open subset of $X$. There is an $n\in\omega$ such that $U_n\subseteq V$, hence $ F_n\subseteq V$.
If $D_1\cap V$ were $G_{\delta}$ in $V$, then $D_1\cap F_n$ were $G_{\delta}$ in $F_n$. As $C_n \subseteq D_1\cap F_n$
it follows that
$D_1\cap F_n$ is dense in $F_n$. Lemma 2.1 then implies that $F_n$ is meager in itself which contradicts the Baire Category
Theorem for $2^\omega$.\par
To finish the proof of the claim it suffices to notice that the countable dense sets $D_0$ and $D_1$ have different (relative)
topological properties in $X$ hence there is no homeomorphism of $X$ sending one to the other, which contradicts the fact
that $X$ is CDH.
\smallskip

Now we are ready to show that $X$ is completely Baire. By Claim 2 and Theorem 1.3, there is a completely metrizable $G\subseteq
X$ which is dense in $X$. Let $D_0$ be any countable dense subset of $G$ (and consequently also a dense subset of $X$.)
Note that $D_0$ has the property that if $E \subseteq D_0$ has no isolated points then $E$ is not $G_\delta$ in $\bar E$, for
if $E$ were $G_\delta$ in $\bar E$ then $E$ would be $G_\delta$ in $\bar E\cap G$, but $\bar E\cap G$ is a $G_\delta$ subset
of $G$, hence, is completely metrizable. However, by Baire Category Theorem this does not happen.\par
Aiming toward a contradiction again, assume that $X$ is not completely Baire. That is, there is a closed set $F\subseteq X$
which is meager in itself. By Lemma 2.1 there is a countable dense $C\subseteq F$ which is $G_\delta$ in $F$. Let $D_1=
C\cup( D_0\setminus F)$. The set $D_1$ is clearly a countable dense subset of $X$ and has the property that
 there is a subset of it without isolated points which is $G_\delta$ in its closure ($C$ being a witness to this.)\par
So, again, the countable dense sets $D_0$ and $D_1$ have different (relative)
topological properties in $X$ hence there is no homeomorphism of $X$ sending one to the other contradicting the
countable dense homogeneity of $X$.\hfill\qed
\enddemo

\bigskip
\proclaim{Corollary 2.4} Every Borel CDH space $X$ is completely metrizable.
\endproclaim

\demo{Proof} Follows directly from Theorem 2.3 and Theorem 1.4.\hfill\qed
\enddemo

\bigskip
\proclaim{Corollary 2.5} Let $X$ be a zero-dimensional Borel CDH space without isolated points. Then  $X$ is homeomorphic
to one of the following five spaces: $2^{\omega}$, $\omega^{\omega}$, $2^{\omega}\smallsetminus\{0\}$,
$\omega^{\omega}\oplus 2^{\omega}$ and $\omega^{\omega}\oplus 2^{\omega}\smallsetminus\{0\} $.
\endproclaim

\demo{Proof} By the previous corollary $X$ is completely metrizable. By Lemma 2.2 $X=L \oplus R$, where
$L$ is locally compact without isolated points and
$R$ has the property that every compact subset or $R$ is nowhere dense in $R$. By Theorem 1.5, 
$R$ is either empty or homeomorphic to $\omega^\omega$ and $L$ (if non-empty) is homeomorphic
either to  $2^{\omega}$ or $2^{\omega}\smallsetminus\{0\}$ depending on whether it is compact
 or not.\hfill\qed
\enddemo

\bigskip
A natural question is whether the above results can be extended beyond analytic or Borel sets.
The answer depends on set theoretic assumptions. For possible extensions note that all arguments 
presented so far use only the validity of Theorems 1.1, 1.3, 1.4
and only countable Axiom of Choice, a consequence of the Axiom of Dependent Choice.

\bigskip
\proclaim{Corollary 2.7} (i) {\rm (PD)} Every projective CDH space is completely metrizable.\newline
(ii) {\rm (AD)} All separable metric CDH spaces are completely metrizable.
\endproclaim

\bigskip

So in particular, it is consistent with ZF that every zero-dimensional metric CDH space without isolated point is 
homeomorphic to one of the following spaces: $2^{\omega}$, $\omega^{\omega}$, $2^{\omega}\smallsetminus\{0\}$,
$\omega^{\omega}\oplus 2^{\omega}$ and $\omega^{\omega}\oplus 2^{\omega}\smallsetminus\{0\} $.

\bigskip

The conclude the section we show that the Theorem 2.3 and
Corollary 2.4 are consistently sharp by proving the following:

\bigskip
\proclaim{Theorem 2.6}{\rm (MA +  $\neg $CH + $\omega_1 =\omega_1^L$)} Let $X$ be an $\aleph_1$-dense subset of $2^\omega$. Then:
(i) $X$  is a co-analytic meager in itself CDH space.\newline
(ii) $2^\omega\setminus X$ is an analytic completely Baire CDH space which is not completely metrizable.
\endproclaim

\demo{Proof} A theorem of Martin and Solovay (see [Mi]) states that, assuming MA +  $\neg $CH + $\omega_1 =\omega_1^L$ , 
every set of reals of size $\aleph_1$ is co-analytic. MA implies that $X$ is meager in itself. It is easy to see
that $2^\omega\setminus X$ is completely Baire and not completely metrizable (and of course analytic). \par
The fact that both $X$ and $2^\omega\setminus X$ are CDH follows directly from Lemma 3.1 of [BB].\hfill\qed
\enddemo

\bigskip

\head 3. Products of CDH spaces \endhead

\bigskip
Theorem 2.6 can be used to see that products of CDH spaces need not be CDH. In fact, if $X$ is a meager in itself
CDH metric space, it is easy to see that $X\times \Bbb R$ is not CDH\footnote{The authors are not aware
of a ZFC example of two metric CDH spaces whose product is not CDH.}. On the other hand, infinite products of
spaces which are not CDH can be CDH, an example being the Hilbert cube $[0,1]^\omega$ ([Fo]). Lawrence [La] showed
that $X^\omega$ is homogeneous, for every $X\subseteq 2^\omega$ (see also [DP]) answering half of Question 388. of [ZH3].
The other half  asks for which $X\subseteq 2^\omega$ is $X^\omega$  CDH. It was known that not for all
as Fitzpatrick and Zhou in [FZ2] showed that $\Bbb Q^\omega$ is not CDH, where $\Bbb Q$ denotes the space of rational
numbers. In this section we characterize those Borel subsets of $2^\omega$ whose power is CDH.

\proclaim{Theorem 3.1} Let $X$ be a separable metric space such that $X^\omega$ is CDH.
Then $X$ is a Baire space.
\endproclaim
 
\demo{Proof} The proof of this theorem is quite analogous to the proof of Claim 2  of Theorem 2.3.  
Suppose that $X$ has at least two elements. It suffices to note,
that (1) if $X$ is not Baire then $X^\omega$ is meager in itself, and (2) Every open subset of $X^\omega$
 contains a copy of $2^\omega$.\hfill
\qed

\proclaim{Theorem 3.2} Let $X\subseteq 2^\omega$ be Borel. Then the following are equivalent:
\roster
\item $X^\omega$ is CDH
\item $X$ is $G_\delta$ in $2^\omega$,
\item $|X^\omega|=1$ or $X^\omega$ is homeomorphic to $2^\omega$ or $X^\omega$ is homeomorphic to $\omega^\omega$.
\endroster
\endproclaim

\demo{Proof} (1) implies (2) by Theorem 2.3 as
$X^\omega$ is Borel if and only if $X$ is and, moreover,
 $X^\omega$ is completely metrizable if and only if $X$ is.\par
To see that (2) implies (3) note that if $X$ is $G_\delta$ in $2^\omega$ then  $X^\omega$ is completely metrizable.
Moreover, if $X$ is zero-dimensional then so is $X^\omega$ and $X^\omega$ does not contain any isolated points.
Suppose that $|X^\omega|>1$.
Now, if $X$ is compact then so is $X^\omega$, hence, $X^\omega$ is homeomorphic to $2^\omega$ by Theorem 1.5 (i). If
$X$ is not compact then all compact subsets of $X^\omega$ are nowhere dense  and $X^\omega$ is homeomorphic to 
$\omega^\omega$ by Theorem 1.5 (ii).\par
(3) implies (1), as both $2^\omega$ and $\omega^\omega$ are CDH.\hfill
\qed

\bigskip
Just like in the previous section this theorem can be strengthen assuming PD or AD.
The following question, however, remains open.

\bigskip
\proclaim{Question 3.2} Is there a non-$G_\delta$ subset of $2^\omega$ such that $X^\omega$ is CDH?
\endproclaim

\bigskip
We will conclude this section and the paper by considering uncountable products.
Recall that a family $\Cal F\subseteq [\omega]^\omega $ is {\it
centered } if every non-empty finite subfamily of $\Cal F$ has an infinite intersection.
An infinite set $A\subseteq\omega$ is a {\it pseudo-intersection } of a family $\Cal
F\subseteq [\omega ]^\omega$  if
$A\setminus F$ is finite for every $F\in \Cal F$. The cardinal invariant $\frak p$ is defined as the minimal cardinality
of a centered family  $\Cal F\subseteq[\omega]^{\omega}$  which has no infinite pseudo-intersection.\par
 
Stepr$\bar{\text{a}}$ns and Zhou in [SZ] showed that $2^\kappa$ is CDH for every $\kappa <\frak p$
and asked whether $2^\frak p$ is provably not CDH. We show that it follows from known results that
the answer is positive.

\bigskip
\proclaim{Theorem 3.3} $\frak p= \min \{\kappa: 2^\kappa$ is not CDH$\}$.
\endproclaim

\demo{Proof} The fact that $\min \{\kappa: 2^\kappa$ is not CDH$\}\leq \frak p$ was proved in [SZ]. In [Ma] and [HS]
it is shown that there is a countable dense set $ D\subseteq 2^\frak p$ and
a point $x\in  2^\frak p$ such that no sequence in $D$ converges to $x$. 
On the other hand, it is easy to construct a countable dense set $C\subseteq 2^\frak p$ such that for every
$c\in C$ there is a sequence $\langle c_n:n\in\omega\rangle \subseteq C\setminus \{c\}$ converging to $c$.\par
Now, notice that there is no homeomorphism of $2^\frak p$ sending $C$ to $D\cap\{x\}$ as if $c=h^{-1}(x)$
and $\langle c_n:n\in\omega\rangle \subseteq C\setminus \{c\}$ a sequence converging to $c$, then
the sequence $\langle h(c_n):n\in\omega\rangle $ does not converge to $x$ contradicting continuity of $h$.
\hfill
\qed

\bigskip
{\bf Acknowledgments.} \rm The work contained in this paper is part of the second author's Master's thesis
at the Universidad Michoacana de San Nicol\'{a}s de Hidalgo, written under the supervision of the first author.
 The first author wishes to thank A. Louveau, I. Farah, J. van Mill and S. Todor\v{c}evi\'{c} for bibliographical 
information and fruitful discussion.

\bigskip

\bigskip

\noindent{Instituto de matem\'{a}ticas, UNAM}
\noindent{Unidad Morelia}
\smallskip
\noindent{A. P. 61-3}
\smallskip
\noindent{Xangari}
\smallskip
\noindent{C. P. 58089, Morelia, Mich., M\'{e}xico}
\smallskip
\noindent{\it michael\@matmor.unam.mx , bzamora\@matmor.unam.mx}
\smallskip

\newpage

\end
\bigskip


\Refs
\widestnumber\key{sh-400a}


\smallskip
\ref \nofrills\key BB
\manyby S.~Baldwin and R~.E.~Beaudoin
\paper  Countable dense homogeneous spaces under {M}artin's axiom 
\jour   Israel J. Math. 
\vol    65
\pages  153--164  
\yr     1989
\endref

\smallskip
\ref \nofrills\key Be
\manyby R~.B.~Bennett
\paper  Countable dense homogeneous spaces 
\jour   Fun. Math 
\vol    74
\pages  189--194 
\yr     1972
\endref

\smallskip
\ref \nofrills\key Br
\manyby L.~E.~J.~Brower
\paper  Some Remarks on the coherence type $\eta$ 
\jour   Proc. Akad. Amsterdam  
\vol    15
\pages  1912 
\yr     1256--1263
\endref

\smallskip
\ref \nofrills\key DP
\manyby A.~Dow and E.~Pearl
\paper  Homogeneity in Powers of zero.dimensional, first-countable spaces 
\jour   Proc. AMS 
\vol    125
\pages  2503--2510
\yr     1997
\endref

\smallskip
\ref \nofrills\key Fi
\manyby  B.~Fitzpatrick Jr.
\paper   A note on countable dense homogeneity
\jour    Fund. Math.
\vol     75
\pages   3--4
\yr      1972
\endref

\smallskip
\ref \nofrills\key FL
\manyby  B.~Fitzpatrick, Jr.and N.~F.~Lauer
\paper   Densely homogeneous spaces. {I}
\jour    Houston J. Math.
\vol     13
\pages   19--25
\yr      1987
\endref

\smallskip
\ref \nofrills\key FZ1
\manyby  B.~Fitzpatrick Jr. and H.-X.~Zhou
\paper   Densely homogeneous spaces II
\jour    Houston J. Math.
\vol     14
\pages   57--68
\yr      1988
\endref

\smallskip
\ref \nofrills\key  FZ2
\manyby  B.~Fitzpatrick Jr. and H.-X.~Zhou
\paper   Countable dense homogeneity and the Baire property
\jour    Topology and its Applications
\vol     43
\pages    1--14
\yr      1992
\endref

\ref\nofrills
\key FZ3
\by B.~Fitzpatrick Jr. and H.-X.~Zhou
\paper Some Open Problems in Densely Homogeneous Spaces
\inbook in  Open problems in Topology (ed. J.~van Mill and M.~Reed) 
\yr 1984
\pages 251--259
\finalinfo North-Holland, Amsterdam 
\endref



\smallskip
\ref \nofrills\key Fo
\manyby  M.~Fort
\paper   Homogeneity of infinite products of manifolds with boundary
\jour    Pacific J. Math
\vol     12
\pages   879--884
\yr      1962
\endref

\smallskip
\ref \nofrills\key Fr
\manyby  M.~Fr\'{e}chet
\paper   Les dimension d'unensemble abstrait
\jour     Math. Ann
\vol     68
\pages   145--168
\yr      1910
\endref

\smallskip
\ref \nofrills\key HS
\manyby  M. Hru\v{s}\'{a}k and J. Steprans
\paper   Cardinal invariants related to sequential separability
\jour    Suri{\-}kai{\-}sekikenkiusho Kokyuroku
\vol     1202
\pages   66--74
\yr      2001
\endref

\smallskip
\ref \nofrills\key Ka
\manyby  A~.Kanamori
\book   The Higher Infinite
\yr     1994
\finalinfo   Springer-Verlag
\endref

\smallskip
\ref \nofrills\key Ke
\manyby  A~.S.~Kechris
\book   Classical Descriptive Set Theory
\yr     1995
\finalinfo   Springer-Verlag
\endref


\smallskip
\ref \nofrills\key KLW
\manyby  A.~S.~Kechris, A.~Louveau and W.~H.~Woodin
\paper   The Structure of $\sigma$-ideals of Compact Sets
\jour    Trans. AMS
\vol     301
\pages   263--288
\yr      1987
\endref


\smallskip
\ref \nofrills\key Ku
\manyby K.~Kunen
\book   Set Theory, An Introduction to Independence Proofs
\yr     1990
\finalinfo   North Holland
\endref

\smallskip
\ref \nofrills\key La
\manyby  B.~Lawrence
\paper   Homogeneity in powers of subspaces of the real line
\jour    Trans. AMS
\vol     350
\pages   3055--3064
\yr      1998
\endref


\smallskip
\ref \nofrills\key Ma
\manyby  M.~V.~Matveev
\paper   Cardinal $\frak p$ and a theorem of Pelczynski
\jour    (preprint)
\vol     
\pages   
\yr      
\endref

\smallskip
\ref \nofrills\key vM1
\manyby J.~van Mill
\paper  Strong local homogeneity does not imply countable dense homogeneity 
\jour   Proc. AMS 
\vol    84
\pages  143--148 
\yr     1982
\endref

\smallskip
\ref \nofrills\key vM2
\manyby J.~van Mill 
\book The Infinite-Dimensional Topology of Function Spaces
\finalinfo  North Holland  
\yr 2001
\endref

\smallskip
\ref \nofrills\key Mi
\manyby A.~W.~Miller 
\book  Descriptive Set Theory and Forcing
\finalinfo  Springer, Lecture Notes in Logic 4
\yr 1995
\endref

\smallskip
\ref \nofrills\key Sa
\manyby W.~L.~Saltsman
\paper  Concerning the existence of a connected, countable dense homogeneous subset of the plane which is not strongly locally homogeneous 
\jour   Topology Proceedings 
\vol    16
\pages  137--176
\yr     1991
\endref

\smallskip
\ref \nofrills\key SZ
\manyby J.~Steprans, H.-X.~Zhou
\paper  Some Results on CDH Spaces
\jour   Topology and its Applications 
\vol    28
\pages  147--154 
\yr     1988
\endref

\smallskip
\ref \nofrills\key Zh
\manyby H.-X.~Zhou
\paper  Two applications of set theory to homogeneity 
\jour   Questions Answers Gen. Topology 
\vol    6
\pages  49--56
\yr     1988
\endref

\endRefs
